\documentclass{amsart}
\usepackage{amssymb}
\usepackage{latexsym}
\usepackage{amscd}
\newtheorem{thm}{Theorem}%[subsection]
\newtheorem{tom}{Theorem}
\newtheorem{lem}{Lemma}%[subsection]
\newtheorem{cor}{Corollary}%[subsection]
\newtheorem{corol}{Corollary}%[subsection]
\newtheorem{prop}{Proposition}%[subsection]
\begin{document}
\renewcommand{\abstractname}{Abstract}
\begin{abstract}
In this note we show the following result using the integral-geometric formula of R. Howard: Consider the totally geodesic $\mathbb{R}P^{2m}$ in $\mathbb{C}P^n$. Then it minimizes volume among the isotropic submanifolds in the same $\mathbb{Z}/2$ homology class in $\mathbb{C}P^n$ (but not among all submanifolds in this $\mathbb{Z}/2$ homology class). Also the totally geodesic $\mathbb{R}P^{2m-1}$ minimizes volume in its Hamiltonian deformation class in $\mathbb{C}P^n$. As a corollary we'll give estimates for volumes of Lagrangian submanifolds in complete intersections in $\mathbb{C}P^n$.
\end{abstract}   
\title[Isotropic submanifolds...]{Volume minimization and estimates for certain isotropic submanifolds in complex projective spaces}     
\author{Edward Goldstein}
\maketitle

\section{Introduction}
On a K\"ahler $n$-fold $M$ there is a class of {\it isotropic} submanifolds. Those are submanifolds of $M$ on which the K\"ahler form $\omega$ of $M$ vanishes. The maximal dimension of such a submanifold is $n$ (the middle dimension) in which case it is called {\it Lagrangian}.\\
In this papers we'll exhibit global volume-minimizing properties among isotropic competitors for certain submanifolds of the complex projective space. In general global volume-minimizing properties of minimal/Hamiltonian stationary Lagrangian/isotropic submanifolds in K\"ahler (particularly K\"ahler-Einstein) manifolds are still poorly understood. In dimesion $2$ there is a result of Schoen-Wolfson \cite{ScW} (extended to isotropic case by Qiu in \cite{Qiu}) which shows existence of Lagrangian cycles minimizing area among Lagrangians in a given homology class. Still it is not clear whether a {\it given} minimal Lagrangian has any global volume-minimizing properties.\\
The only instance where we have a clear cut answer to global volume-minimizing problem is Special Lagrangian submanifolds which are homologically volume-mimizing in Calabi-Yau manifolds \cite{HaL}. In K\"ahler-Einstein manifolds of negative scalar curvature, besides geodesics on Riemann surfaces of negative curvature, we have some examples \cite{Lee} of minimal Lagrangian submanifolds which are homotopically volume-minimizing. The author has a program for studying homotopy volume-minimizing properties for Lagrangians in K\"ahler-Einstein manifolds of negative scalar curvature \cite{Gold1}, but so far there are no satisfactory results.\\
In positive curvature case there is a result of Givental-Kleiner-Oh which states that the canonical totally geodesic $\mathbb{R}P^{n}$ in $\mathbb{C}P^{n}$ minimizes volume in its Hamiltonian deformation class, \cite{Giv}. The proof uses integral geometry and Floer homology to study intersections for Hamiltonian deformations of $\mathbb{R}P^{n}$. Those arguments can be generalized to products of Lagrangians in a product of symmetric K\"ahler manifolds, \cite{IOS}. There is a related conjecture due to Oh that the Clifford torus minimizes volume in its Hamiltonian deformation class in $\mathbb{C}P^{n}$, \cite{Oh}. Some progress towards this was obtained in \cite{Gold2}. Also general lower bounds for volumes of Lagrangians in a given Hamiltonian deformation class in $\mathbb{C}^n$ were obtained in \cite{Vit}.\\
In this note we extend and improve the result of Givental-Kleiner-Oh to isotropic totally geodesic $\mathbb{R}P^k$ sitting canonically in $\mathbb{C}P^n$. Our main result is the following theorem:
\begin{tom}
Consider the totally geodesic $\mathbb{R}P^{2m}$ in $\mathbb{C}P^n$. Then it minimizes volume among the isotropic submanifolds in the same $\mathbb{Z}/2$ homology class in $\mathbb{C}P^n$ (but not among all submanifolds in this $\mathbb{Z}/2$ homology class). Also consider the totally geodesic $\mathbb{R}P^{2m-1}$ in $\mathbb{C}P^n$. Then it minimizes volume in its Hamiltonian deformation class.
\end{tom}
A corollary of this is:
\begin{corol}
Let $f_1,\ldots,f_k$ be real homogeneous polynomials of odd degree in $n+1$ variables with $2m+k=n$. Let $N$ be the zero locus of $f_i$ in $\mathbb{C}P^n$ and $L$ be their real locus. Then $vol(L) \leq \Pi deg(f_i) vol(\mathbb{R}P^{2m})$ and if $L'$ is a Lagrangian submanifold of $N$ homologous mod $2$ to $L$ in $N$ then $vol(L') \geq vol(\mathbb{R}P^{2m})$.
\end{corol}
\section{A formula from integral geometry}
In this section we establish a formula from integral geometry for volumes of isotropic submanifolds of $\mathbb{C}P^n$ following the exposition in R. Howard \cite{How}.\\
In our case the group $SU(n+1)$ acts on $\mathbb{C}P^n$ with a stabilizer $K \simeq U(n)$. Thus we view $\mathbb{C}P^{n}=SU(n+1)/K$ and the Fubini-Study metric is induced from the bi-invariant metric on $SU(n+1)$. Let $P^{2m}$ be an isotropic submanifold of $\mathbb{C}P^n$ of dimension $2m$ and let $Q$ be a linear $\mathbb{C}P^{n-m} \subset \mathbb{C}P^n$. For a point $p \in P$ and $q \in Q$ we define an angle $\sigma(p,q)$ between the tangent planes $T_pP$ and $T_qQ$ as follows: First we choose some elements $g$ and $h$ in $SU(n+1)$ which move $p$ and $q$ respectively to the same point $r \in \mathbb{C}P^n$. Now the tangent planes $g_{\ast}T_pP$ and $h_{\ast}T_qQ$ are in the same tangent space $T_r \mathbb{C}P^n$ and we can define an angle between them as follows: take an orthonormal basis $u_1\ldots u_{2m}$ for $g_{\ast}T_pP$ and an orthonormal basis $v_1 \ldots v_{2n-2m}$ for $h_{\ast}T_qQ$ and define \[\sigma(g_{\ast}T_pP,h_{\ast}T_qQ)= |u_1 \wedge \ldots \wedge v_{2n-2m}|\]
The later quantity $\sigma(g_{\ast}T_pP,h_{\ast}T_qQ)$ depends on the choices $g$ and $h$ we made. To mend this we'll need to average this out by the stabilizer group $K$ of the point $r$. Thus we define:
\[\sigma(p,q)=\int_{K} \sigma(g_{\ast}T_pP, k_{\ast}h_{\ast}T_qQ) dk\]
Since $SU(n+1)$ acts transitively on the Grassmanian of isotropic planes and the complex planes in $\mathbb{C}P^n$ we conclude that this angle is a constant depending just on $m$ and $n$:
\[\sigma(p,q)=C_{m,n}\]
There is a following general formula due to R. Howard \cite{How}:
\[ \int_{SU(n+1)} \#(P \bigcap gQ) dg= \int_{P \times Q} \sigma(p,q) dp dq= C_{m,n} vol(P)vol(Q)\]
Here $\#(P \bigcap gQ)$ is the number of intersection points of $P$ with $gQ$, which is finite for a generic $g \in SU(n+1)$. To use the formula we need to have some control over the intersection pattern of $P$ and $gQ$. We have the following lemma:
\begin{lem}
\label{tech}
Let $P$ be the totally geodesic $\mathbb{R}P^{2m} \subset \mathbb{C}P^n$, let $Q=  \mathbb{C}P^{n-m} \subset \mathbb{C}P^n$. Let $g \in SU(n+1)$ s.t. $P$ and $gQ$ intersect transversally. Then $ \#(P \bigcap gQ)=1$. Also let $f_1,\ldots,f_k$ be real homogeneous polynomials in $n+1$ variables with $2m+k=n$ and let $P'$ be their real locus. If $P'$ is transversal to $gQ$ then $ \#(P' \bigcap gQ) \leq \Pi deg(f_i)$.
\end{lem}
{\bf Proof:} For the first claim we have $gQ$ is given by an $(n-m+1)$-plane $H \subset \mathbb{C}^{n+1}$ and hence it is a zero locus of $m$ linear equations on $\mathbb{C}^{n+1}$. Hence $(P \bigcap gQ)$ is cut out by $2m$ linear equations in $\mathbb{R}P^{2m}$.\\
For the second claim we note that as before $gQ \bigcap \mathbb{R}P^n$ is the zero locus of $2m$ linear polymonials $h_1,\ldots, h_{2m}$ on $\mathbb{R}P^n$. Moreover $P'$ is a zero locus of $f_1,\ldots,f_{n-2m}$ on $\mathbb{R}P^n$. For generic $g \in SU(n+1)$ we'll have that $gQ$ and $P'$ intersect transversally in $\mathbb{R}P^n$. By Bezout's theorem (see \cite{GH}, p. 670) the common zero locus of $h_1,\ldots, h_{2m}$ and $f_1,\ldots,f_{n-2m}$ is $\mathbb{C}P^n$ is $\Pi deg(f_i)$ points. Now $P' \bigcap gQ$ is a part of this locus, hence $ \#(P' \bigcap gQ) \leq \Pi deg(f_i)$.
 
\section{Proof of the volume minimization}

Now we can prove the result stated in the Introduction: 
\begin{thm}
Consider the totally geodesic $\mathbb{R}P^{2m}$ in $\mathbb{C}P^n$. Then it minimizes volume among the isotropic submanifolds in the same $\mathbb{Z}/2$ homology class in $\mathbb{C}P^n$ (but not among all submanifolds in this $\mathbb{Z}/2$ homology class). Also consider the totally geodesic $\mathbb{R}P^{2m-1}$ in $\mathbb{C}P^n$. Then it minimizes volume in its Hamiltonian deformation class.
\end{thm}
{\bf Proof:} Let $P$ be an isotropic submanifold homologous to $\mathbb{R}P^{2m}$ mod $2$ and let $Q=\mathbb{C}P^{n-m}$. By Lemma \ref{tech} the intersection number mod $2$ of $P$ and $gQ$ is $1$. Hence the formula in the previous section tells that 
\[C_{m,n} vol(P)vol(Q)= \int_{SU(n+1)} \#(P \bigcap gQ) dg \geq vol(SU(n+1))\]
and
\[C_{m,n} vol(\mathbb{R}P^{2m})vol(Q)= \int_{SU(n+1)} \#(\mathbb{R}P^{2m} \bigcap gQ) dg = vol(SU(n+1))\]
and this proves the first part. We also note that that $\mathbb{C}P^1$ is homologous to $\mathbb{R}P^2$ mod $2$ in $\mathbb{C}P^n$ but \[vol(\mathbb{C}P^1) <vol(\mathbb{R}P^2) \] 
The second assertion will follow from the first one. Consider $\mathbb{C}^{n+1}$ and a unit sphere $S^{2n+1} \subset \mathbb{C}^{n+1}$. We have a natural circle action on $S^{2n+1}$ (multiplication by unit complex numbers). Let the vector field $u$ be the generator of this action. We have a $1$-form $\alpha$ on $S^{2n+1}$, \[ \alpha(v) = u \cdot v \] Also $d\alpha=2\omega$ where $\omega$ is the K\"ahler form of $\mathbb{C}^{n+1}$. The kernel of $\alpha$ is the {\it horizontal distribution}. We have a Hopf map $\rho: S^{2n+1} \mapsto \mathbb{C}P^n$. We have $\mathbb{R}P^{2m-1} \subset \mathbb{C}P^n$ and $S^{2m-1} \subset  S^{2n+1}$ which is a horizontal double cover of  $\mathbb{R}P^{2m-1}$. \\
Let $f$ be a (time-dependent) Hamiltonian function on  $\mathbb{C}P^n$. Then we can lift it to a Hamiltonian function on $\mathbb{C}^{n+1}-(0)$ and its Hamiltonian vector field $H_f$ is horizontal on $S^{2n+1}$. Consider now the vector field \[w=-2f \cdot u+H_f\]
The vector field $w$ is $S^1$-invariant. We also have:
\begin{prop}
The Lie derivative $L_{w}\alpha= 0$
\end{prop}
{\bf Proof:} We have \[L_{w}\alpha= d(i_w \alpha)+i_w d \alpha= -2df+2df\]
Let now $\Phi_t$ be the time $t$ flow of $w$ on $S^{2n+1}$ and let $\Xi_t$ be the Hamiltonian flow of $f$ on $\mathbb{C}P^n$. Then $\Phi_t(S^{2m-1})$ is horizontal and isotropic and it is a double cover of $\Xi_t(\mathbb{R}P^{2m-1})$. Hence 
\[vol(\Phi_t(S^{2m-1}))=2vol(\Xi_t(\mathbb{R}P^{2m-1})) \]
Let $S_t=\Phi_t(S^{2m-1})$. We build a suspension $\Sigma S_t$ of $S_t$ in $S^{2n+3} \subset \mathbb{C}^{n+2}$, \[\Sigma S_t = \big((\sin \theta \cdot x, \cos \theta) \in  \mathbb{C}^{n+2}=  \mathbb{C}^{n+1} \oplus  \mathbb{C}| 0 \leq \theta \leq \pi ~ , ~  x \in S_t \big)\]
One immediately verifies that $\Sigma S_t$ is horizontal and it is a double cover of an isotropic submanifold $L_t$ (with a conical singularity) of $\mathbb{C}P^{n+1}$ with $L_0=\mathbb{R}P^{2m}$. Also one readily checks that \[vol(\Sigma S_t)= vol(S_t) \cdot \int_{\theta=0}^{\pi} \sin^{2m-1} \theta \ d \theta\] 
Hence
\[2 vol(L_t)=vol(\Sigma S_t)=  2vol(\Xi_t(\mathbb{R}P^{2m-1})) \cdot \int_{\theta=0}^{\pi} \sin^{2m-1} \theta \ d \theta \]
Now the first part of our theorem implies that $vol(L_t) \geq vol(L_0)$. Hence we conclude that $vol(\Xi_t(\mathbb{R}P^{2m-1})) \geq vol(\mathbb{R}P^{2m-1})$.
 Q.E.D. \\
{\bf Remark:} One notes from the proof that for $\mathbb{R}P^{2m-1}$ it would be suffient to use exact deformations by isotropic immersions of 
$\mathbb{R}P^{2m-1}$. A family $L_t$ of isotropic immersions of 
$\mathbb{R}P^{2m-1}$ is called {\it exact} if the $1$-form $i_v \omega$ is exact when restricted to each element of the family. Here $v$ is the deformation vector field and $\omega$ is the symplectic form. Thus embeddedness is not important for the conclusion of the theorem.\\
The theorem has the following corollary:
\begin{cor}
Let $f_1,\ldots,f_k$ be real homogeneous polynomials of odd degree in $n+1$ variables with $2m+k=n$. Let $N$ be the zero locus of $f_i$ in $\mathbb{C}P^n$ and $L$ be their real locus. Then $vol(L) \leq \Pi deg(f_i) vol(\mathbb{R}P^{2m})$ and if $L'$ is a Lagrangian submanifold of $N$ homologous mod $2$ to $L$ in $N$ then $vol(L') \geq vol(\mathbb{R}P^{2m})$.
\end{cor}
{\bf Proof:} We note that $N$ is a complex $2m$-fold and $L$ is its Lagrangian submanifold. Since the degrees of $f_i$ are odd, we have by adjunction formula that $L$ and $\mathbb{R}P^{2m}$ represent the same homology class in $H_{2m}(\mathbb{R}P^n,\mathbb{Z}/2)$. Let $Q$ be a linear $\mathbb{C}P^{n-m}$ in $\mathbb{C}P^n$ and $g \in SU(n+1)$. The intersection munber mod $2$ of $gQ$ with $L'$ is $1$. We have that 
\[C_{m,n} vol(\mathbb{R}P^{2m})vol(Q)= \int_{SU(n+1)} 1 dg\]
\[C_{m,n} vol(L')vol(Q)= \int_{SU(n+1)} \#(L' \bigcap gQ) dg \]
Also using Lemma \ref{tech}:
\[C_{m,n} vol(L)vol(Q)= \int_{SU(n+1)} \#(L \bigcap gQ) dg \leq \Pi deg(f_i) vol(SU(n+1))\] and our claims follow.  Q.E.D. \\

egold@ias.edu


\begin{thebibliography}{99}


\bibitem[Giv]{Giv} A. Givental: The Nonlinear Maslov index, London Mathematical Society Lecture Note Series 15 (1990), 35-43

\bibitem[Gold1]{Gold1} Edward Goldstein: Strict volume-minimizing properties for Lagrangian submanifolds in complex manifolds with positive canonical bundle, math.DG/0301191

\bibitem[Gold2]{Gold2} Edward Goldstein: Some estimates related to Oh's conjecture for the Clifford tori in $\mathbb{C}P^n$, math.DG/0311460

\bibitem[GH]{GH}P. Griffiths, J. Harris, ``Principles of Algebraic geometry,'' Wiley and Sons, 1978.

\bibitem[HaL]{HaL} R. Harvey and H. B. Lawson : Calibrated Geometries, Acta
Math. 148, 47-157 (1982). 

\bibitem[How]{How} Howard, Ralph: The kinematic formula in Riemannian homogeneous spaces. Mem. Amer. Math. Soc. 106 (1993), no. 509, vi+69 pp.

\bibitem[IOS]{IOS} Hiroshi Iriyeh, Hajime Ono, Takashi Sakai: Integral Geometry and Hamiltonian volume minimizing property of a totally geodesic Lagrangian
    torus in $S^2 \times S^2$, Proc. Japan Acad. Ser. A Math. Sci. 79 (2003), no. 10, 167-170

\bibitem[Lee]{Lee} Y.-I. Lee: Lagrangian minimal surfaces in K\"ahler-Einstein surfaces of negative scalar curvature. Comm. Anal. Geom. 2 (1994), no. 4, 579--592.

\bibitem[Oh]{Oh} Y.-G.~Oh: Mean curvature vector and
  symplectic topology of Lagrangian submanifolds in Einstein-K\"ahler
  manifolds, Math.~Z.~{\bf 216}, 471-482 (1994). 

\bibitem[Qiu]{Qiu} Qiu, Weiyang: Interior regularity of solutions to the isotropically constrained Plateau problem. Comm. Anal. Geom. 11 (2003), no. 5, 945--986

\bibitem[ScW]{ScW} Schoen, R.; Wolfson, J.: Minimizing area among Lagrangian surfaces: the mapping problem. J. Differential Geom. 58 (2001), no. 1, 1--86.

\bibitem[Vit]{Vit} C. Viterbo: Metric and isoperimetric problems in symplectic
   geometry. J. Amer. Math. Soc. 13 (2000), no. 2, 411--431 

\end{thebibliography}
\end{document}